\documentclass{article}
\usepackage{amssymb}
\usepackage{amsmath}
\usepackage{graphicx}
\usepackage{epstopdf}
\usepackage{appendix}

\begin{document}
\setlength\parindent{0pt}

\begin{titlepage}
\begin{center}

\textsc{\LARGE American University in Bulgaria}\\[4.5cm]
\textsc{\Large Senior Thesis}\\[0.7cm]

{ \huge \bfseries Tessellations}\\[1.0cm]

\begin{minipage}{0.4\textwidth}
\begin{flushleft} \large
\emph{Author:}\\
Plamen \textsc{Dimitrov}
\end{flushleft}
\end{minipage}
\begin{minipage}{0.4\textwidth}
\begin{flushright} \large
\emph{Advisor:} \\
Prof.~Orlin \textsc{Stoytchev}
\end{flushright}
\end{minipage}

\vfill
{\large \today}

\end{center}
\end{titlepage}

\vspace*{2.5cm}
\begin{abstract}
\noindent This work presents the tessellations and polytopes from the perspective of both n-dimensional geometry and abstract symmetry groups. It starts with a brief introduction to the terminology and a short motivation. In the first part, it engages in the construction of all regular tessellations and polytopes of n dimensions and extends this to the study of their quasi-regular and uniform generalizations. In the second part, the symmetries of polytopes and tessellations are considered and the Coxeter groups and their associated root systems are introduced and classified. In the last part, the algorithms developed for this work are described and their results discussed.
\end{abstract}
\newpage

\tableofcontents
\newpage

\section{Introduction}

\subsection{Motivation and tessellation definitions}

The idea for this work started from the visualization of the sequence of 4-cube cross-sections with a 3-flat (3-dimensional affine space) which results in an octahedron at the center of the 4-cube and in tetrahedrons of different size and truncation magnitude along any of the 4-fold diagonals of the 4-cube. \newline

\begin{figure}[h!]
  \centering
    \includegraphics[width=0.7\textwidth]{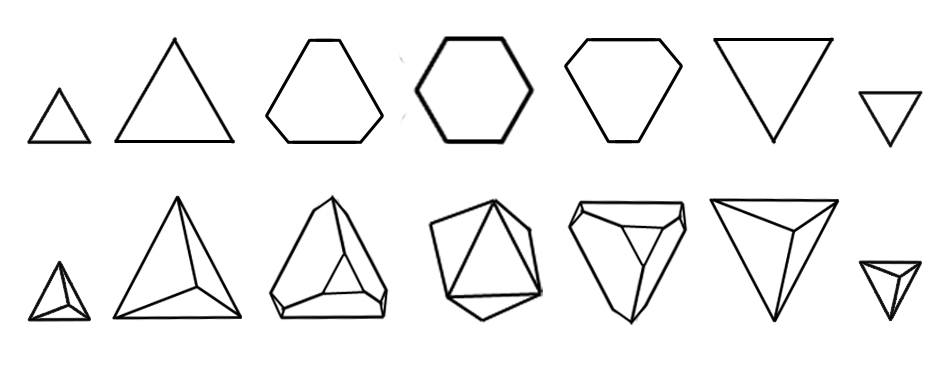}
  \caption{3-dimensional cross-sections along 4-fold symmetry axis.}
\end{figure}

The fact that the three regular polytopes in all dimensions could be obtained from a cross-section with a hypercube and that various uniform polytopes can further be obtained from these together with the fact that the hypercubic tessellation is the only regular tessellation in all dimensions led to the idea to develop an algorithm that could help us identify specific types of symmetric cross-sections with the hypercubic tessellation. \newline

Since there were many different naming conventions used by mathematicians regarding some of the terms studied in this text, this first subsection is devoted to a clarification regarding those. Most mathematicians distinguish among tilings, honeycombs and tessellations, using the last as a generalization of $n$-dimensional case and the first two to thus refer to tessellations on the plane and on a higher dimensional space. Here we employ the same meaning for these terms but avoid their interchangeability which is possible in some sources. \newline

The definition of a tessellation used here is due to Coxeter \cite{coxeter63}. A tessellation is an infinite set of polytopes fitting together to fill an $n$-dimensional space without overlapping so that the $n$-1-facet of each polytope belongs to one and only one other polytope. \newline

Tessellations in $n-1$ dimensions can be considered as $n$-dimensional apeirotopes, i.e. polytopes with infinitely many cells. This will be justified in the next section where the connection between polytopes and tessellations will be examined and it will be shown that a polytope may also be regarded as a tessellation of a given manifold. Before this, a few definitions regarding polytopes are necessary. \newline

\subsection{Polytope definitions}

A polytope $P$ is a geometrical figure bounded by finitely many hyperplanes. It is the $n$-dimensional generalization of polygon (2-polytope), polyhedron (3-polytope), polychoron (4-polytope), etc. with the following properties:

\begin{itemize}
  \item The facets of $P$ are polytopes.
  \item $P$ has facets of every dimension $0,1,...,n-1$ and for brevity a $k$-dimensional facet is called a $k$-facet of $P$.
  \item Every facet of a facet of $P$ is also a facet of $P$.
  \item A convention assumed in this text is the existence of one $n$-facet of $P$ which is $P$ itself.
  \item For every two facets $F_1$ and $F_2$ of $P$, $F_1 \cap F_2$ is also a facet of $P$.
  \item For every two facets $F_1$ and $F_2$ of $P$ there exists a uniquely defined facet $F_1 \vee F_2$ of $P$, namely the smallest facet in terms of inclusion which contains both $F_1$ and $F_2$ .
\end{itemize}

The last two properties originate from the following observation. For an $n$-polytope, each $n$-1-facet lies in a bounding hyperplane. Two or more adjacent bounding hyperplanes intersect at $n$-2-flat and $n$-2-facet of $P$ lies in this $n$-2-flat, three or more adjacent $n$-2-flats intersect at $n$-3-flat and $n$-3-facet lies in this $n$-3-flat and so on as $n$ or more 1-flats or lines intersect in points which are the vertices of the polytope. \newline

There are different naming conventions of the various facets of a polytope besides the obvious for vertex, edge, and face. Examples are cell for 3-facet, hypercell for 4-facet, facet for $n$-1-facet, ridge for $n$-2-facet, and peak for $n$-3-facet. For the sake of any further generalization, only the first part of these will be used. Each constituent which is also a polytope of a lower dimension is additionally referred to as $k$-face, $k$-boundary or similar alternative names but here $k$-facet ($k<n$) will be used for this purpose. \newline

Let $N_r$ denote the number of $r$-facets of a polytope and $N_{pq}$ denote the number of $p$-facets that lie in a $q$-facet if $p<q$ or the number of $p$-facets that pass through a $q$-facet if $p>q$. Then $N_r=N_{rr}=N_{rn}$ and it can be noticed that these numbers are not independent, i.e. $N_{pq}N_{q}=N_{qp}N_{p}$. They will be called configurational numbers and will be used later in this work. \newline

Polytopes may be divided into: convex and non-convex (also known as star-polytopes); regular, quasi-regular and semi-regular; uniform and non-uniform. A set $S \subset \mathbb{R}^n$ is convex if it has the property that for any pair of points $x, y \in S$, the line segment

\[{\lambda x + (1 - \lambda) y : 0 \leq \lambda \leq 1}\]

with end points $x$ and $y$, lies entirely in $S$. For any set $S$, the smallest convex set containing $S$ (the intersection of the family of all convex sets that contain $S$) is called the convex hull of $S$, and is denoted by $conv(S)$. A convex polytope $P$ is defined to be the convex hull of any finite set of points in $\mathbb{R}^n$ \cite{gruenbaum69}. If $F_1^0$, ..., $F_r^0$ are the end-points of the edges of $P$ that meet at $F^0$, then $conv(F_1^0$, ..., $F_r^0)$ is an $n$-1-polytope called the vertex figure of $P$. An $n$-polytope $P$ is regular if its facets are regular and its vertex figures are regular (this definition will be revisited in the next section). Finally, an $n$-polytope $P$ is uniform if its facets are uniform and its vertex figures are all of the same kind ($P$ is vertex-transitive). \newline

The dual of a polytope is defined as a polytope whose $r$-facets correspond to the $(n-r-1)$-facets of the original so that its $n$-1-facets are the vertex figures of the original. The configurational numbers of the $n$-1-facet of the dual polytope then naturally correspond to the configurational numbers of the vertex figure of the original polytope. \newline

Star-polytopes will not be considered in this work.

\section{The regular tessellations and polytopes of n dimensions}

\subsection{The regular tilings and polyhedra}

In the following section we will enumerate the regular polytopes in $n$ dimensions and will use the same principle to define all regular tessellations of $n$-1 dimensions as infinite regular $n$-polytopes. The construction of these is based on Sommerville \cite{sommerville63}. \newline

Regular polyhedra satisfy two conditions: (1) their faces are regular polygons of the same kind; (2) their solid angles are equal. In the classical definition however, the second condition is restated in a stronger form: the vertex figures are regular polygons, which allows for the first condition to be weakened to: their faces are regular polygons. A more modern definition requires an even stronger condition - the polyhedron must be transitive on its flags, i.e. face-, edge-, and vertex-transitive, i.e. all faces, edges and vertices must be the same. The reason for this new definition is that while the requirement for regular faces and vertex figures is sufficient to derive this condition in the usual case, it is no longer valid for the newly introduced abstract polytopes. \newline

The two final conditions to be used for the following construction are weaker: (1) each face has the same number $n$ of edges and vertices; (2) there is the same number $p$ of edges and faces around each vertex. These two numbers $n$ and $p$ will show to be sufficient for the construction of the same final polytopes as with the stronger conditions. From this definition then it follows that $N_{02}=N_{12}=n$ and $N_{10}=N_{20}=p$ where $N_{ij}$ is the number of i-faces incident to a j-face. Since each edge is surrounded by two vertices and faces, $N_{01}=N_{21}=2$. Therefore using $N_{ij}N_{j}=N_{ji}N_{i}$ we arrive at $n N_2 = 2 N_1 = p N_0$. We have to now add one more restriction that will ultimately limit the possible values of $n$ and $p$. \newline

This restriction is Euler's polyhedron formula

\[N_2 - N_1 + N_0 = \chi\]

where $\chi$ denotes a topological invariant called the Euler characteristic. A regular polyhedron then has the Euler characteristic of the sphere $\chi = 2$. To verify this, we need some further tools. The Gauss - Bonnet formula

\[\int_{M}KdA=2\pi \chi(M)\]

relates the Gaussian curvature $K$ to the Euler characteristic $\chi$ as $M$ is any orientable closed surface. Descartes' law of closure defect then states that if the polyhedron is homeomorphic to a sphere (hence not necessarily convex), the total angle defect as the sum of the defects of all the vertices is $4\pi$ which is the total Gaussian curvature of the sphere. The Gauss - Bonnet formula then gives

\[4\pi = \int_{M}KdA=2\pi \chi(M) \Rightarrow \chi(M)=2.\]

The surface of a polyhedron can be considered as a limit case of differentiable surface where the total Gaussian curvature $K$ remains zero at the faces and edges of the polyhedron and is therefore concentrated on the discrete vertex points. Although $K$ is not defined at these points, its integral remains finite and the integral of the total curvature can be replaced with the sum of the angle defect at all vertices. Thus, the number of vertices can be easily derived by dividing $2\pi \chi(M)=4\pi$ with the angle defect at a vertex.

\[N_0(2+\frac{2n}{p}-n)\pi=4\pi\]

The options for the selection of $p$ polygons around a vertex and $n$ vertices around a polygon to get an appropriate angle defect are finally limited to this divisibility criterion. Although this provides with solutions for $p$ and $n$, a more geometric approach which can be applied to a higher dimensional case is considered here.\newline

From Euler's formula for a regular polyhedron $N_2 - N_1 + N_0 = 2$ we calculate

\[N_0=2-N_2+N_1=2-\frac{pN_0}{n}+\frac{pN_0}{2} \Rightarrow N_0(1 +\frac{p}{n}-\frac{p}{2})=2 \Rightarrow N_0=\frac{4n}{2(n+ p)-np}\]

and analogously for the rest we get

\[N_1=\frac{2np}{2(n+p)-np}\mbox{, }N_2=\frac{4p}{2(n+p)-np}\mbox{, }\lambda = \frac{2}{2(n+p)-np}\]

which simplifies to $N_0=2n\lambda, N_1=np\lambda\mbox{, and }N_2=2p\lambda$. \newline

Restricting this to finite regular polyhedra requires finite and positive $\lambda$ which imposes the additional conditions $p>2, n>2$ and $2(n+p)-np>0 \Rightarrow p < \frac{2n}{n-2}$. Therefore, for $n=3 \Rightarrow 2<p<6, n=4 \Rightarrow 2<p<4, n=5 \Rightarrow 2<p<4$. $n \geq 6 $ would imply $2<p<3$ which is not possible, therefore this depletes all options. As a result, all possible regular polyhedrons are

\begin{center}
\begin{tabular}{ r l l l l l }
  n & 3 & 3 & 3 & 4 & 5 \\
  p & 3 & 4 & 5 & 3 & 3 \\
\end{tabular}
\end{center}

These are exactly the Platonic solids. The tetrahedron \{3,3\} is self-dual, while the hexahedron \{4,3\} has the octahedron \{3,4\} as its dual and the icosahedron \{3,5\} is the dual of the dodecahedron \{5,3\}. In a similar way we can take an infinite regular polyhedron by setting $2(n+p)-np=0$ for which $\lambda$ becomes infinite. Since $p>2$ and $n>2$, this leads to 

\begin{center}
\begin{tabular}{ r l l l }
  n & 3 & 4 & 6 \\
  p & 6 & 4 & 3 \\
\end{tabular}
\end{center}

These are the three regular tilings, respectively the triangular \{3,3\}, the hexagonal \{6,3\}, and the square \{4,4\}. \newline

The case $2(n+p)-np<0$ is a third case but in order to understand it better, we must interpret these results first. In the first case, the polyhedron's configuration of vertices, edges and faces is isomorphic to a tessellation of a sphere where all its vertices lie on the surface and its edges are represented by geodesic arcs. To see this, let $S$ denote the total sum of angles of an $n$-gon on the sphere with radius $r$. Since the total angle deficiency for a sphere is

\[\alpha = \int_{M}KdA=\frac{1}{r^2}A(M) \Rightarrow A(M)=r^2\alpha\]

and for an $n$-gon

\[2\pi-\alpha+S=n\pi\Rightarrow \alpha=(S-(n-2)\pi)\]

it follows that the area of the $n$-gon must be $A(M)=(S-(n-2)\pi)r^2$ and since $\sum_{i=1}^{N_2} S_i = \sum_{i=1}^{N_2} S = 2\pi N_0$, and $n N_2 = 2 N_1 = p N_0$, the total area must be

\[((2N_0-(n-2)N_2)\pi r^2 = (2N_0+2 \frac{p N_0}{n}-pN_0)\pi r^2 = \frac{2(n+p)-np}{n}N_0\pi r^2.\]

Then for the case $2(n+p)-np=0$ either $r$, $N_0$ or both have to be infinite in order to have a finite area of the $n$-gon. If $r$ is infinite the sphere becomes a plane else $N_0$ must be infinite which makes all vertices indefinitely close to each other on a finite sphere thus again making Euclidean geometry applicable (with sum of angles of $n$-gon now $S = (n-2)\pi$). The third case follows from the first two. \newline

The last condition $2(n+p)-np<0$ then implies $r^2<0$ where a purely imaginary radius would require trigonometric formulae which hold on a surface with constant negative Gaussian curvature $K$ called Lobachevski sphere or a hyperbolic plane. These three cases then help us conclude that regular polyhedrons are equivalent to tessellations of elliptic, Euclidean and hyperbolic plane. \newline

The above condition $2(n+p)-np>0$ can also be rewritten as

\[2(n+p)-np>0 \Rightarrow \frac{2n}{p}+(2-n)>0 \Rightarrow \frac{2}{p}>1-\frac{2}{n} \Rightarrow \frac{1}{p}+\frac{1}{n}>\frac{1}{2}\]

which is another inequality related to the vertex figure's angle defect.

\subsection{The regular cubic honeycomb and polychora}

In order to continue this construction in n dimensions, we now use the more restrictive classical definition for $n$-dimensional regular polytope with the two conditions: (1) all $n$-1-facets of the polytope must be regular polytopes; (2) all vertex figures of the polytope must be regular polytopes. Although it is possible to consider polytopes with $n$-1-facets that are tessellations of $n$-2-dimensional space, a restriction only to finite polytopes as $n$-1-facets and vertex figures is considered for the purpose of this construction. \newline

To continue the approach from the previous subsection, extend the configurational numbers $N_{pq}$ of the polytope with configurational numbers of the $n$-1-facet $F_{pq}$ and of the vertex figure $V_{pq}$. Then the following equalities hold because of the definition:

\begin{center}
\begin{tabular}{ l l l l }
  $N_{pq} = F_{pq}$ & $(p < q < 3)$, & $N_{p3} = F_p$ & $(p = 0, 1, 2)$ \\
  $N_{pq} = V_{p-1,q-1}$ & $(p > q > 0),$ & $N_{p0} = V_{p-1}$ & $(p = 1, 2, 3)$ \\
\end{tabular}
\end{center}

The above equalities are true for the specific values of $p$. The first equality constitute the relations between the facet and the vertex figure of a 4-polytope. The second is derived from $N_{p3} = F_{p3}=F_p$ since the n-1-facets of the 3-polytope are 3-facets. The third equality makes use of duals $N'$: $V_{pq}=F'_{pq}=N'_{pq}=N_{p-1,q-1}$. Finally, the last equality is similar to the second in terms of the dual polytope. \newline

Furthermore, the number of lines or planes through a point $F_{10}, F_{20}$ in a 3-facet (cell) must be the same as the number of lines or planes $V_{02}, V_{12}$ through each vertex, hence $F_{10}=F_{20}=V_{02}=V_{12}$. This is one of the three numbers $p$, $q$, $r$ which can finally be extracted from

\[F_{02}=F_{12}=N_{02}=N_{12}=p\]
\[F_{10}=F_{20}=V_{02}=V_{12}=q\]
\[V_{10}=V_{20}=N_{21}=N_{31}=r\]

where $p,q$ describe the 3-facet and $q,r$ the vertex figure. They represent the number of vertices or edges of a polygon ($p$), the number of edges or planes through each vertex of a polyhedron ($q$), and the number of 3-facets or cells through each edge ($r$). \newline

As mentioned before, by superimposing all available choices for finite polyhedrons, the following 4-dimensional cases can be constructed:

\begin{center}
\begin{tabular}{ l l l l l l }
  333 & 334 & 335 & 343 & 353 & 433 \\
  434 & 435 & 533 & 534 & 535 &  \\
\end{tabular}
\end{center}

However, they still need to be distinguished in terms of metric (elliptic, Euclidean, and hyperbolic). To do this, we have to examine the dihedral angles of the constituent polyhedrons. The geometric prescription makes use of spherical trigonometry. \newline

\begin{figure}[h!]
  \centering
	$\begin{array}{cc}
		\includegraphics[scale=0.25]{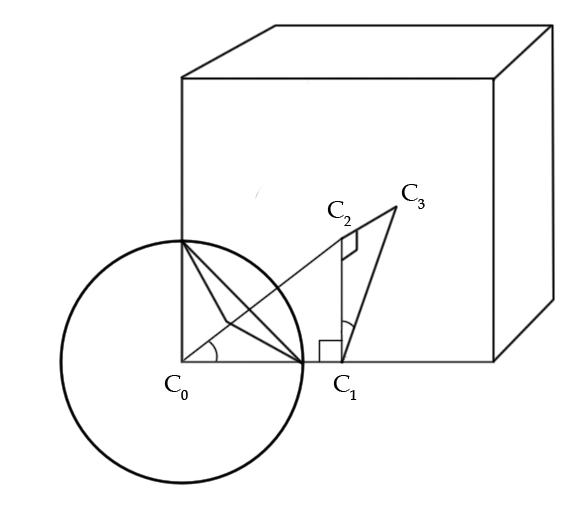} &
		\includegraphics[scale=0.3]{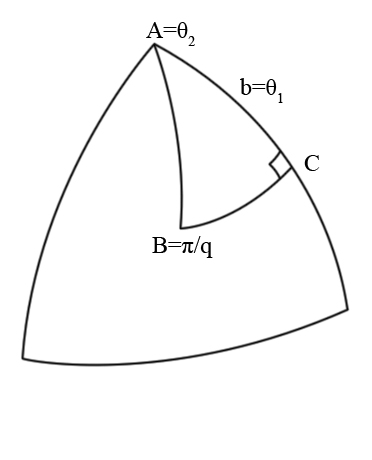}
	\end{array}$
  \caption{Dihedral angles of a polyhedron.}
\end{figure}

Let $C_0$ be any vertex a regular polyhedron p,q. Let $C_1$ be the mid-point of an edge through $C_0$, $C_2$ the center of a face though $C_0C_1$, and $C_3$ the center of the polyhedron. Denote $\angle C_1C_0C_2 = \theta_1$ and $\angle C_2C_1C_3 = \theta_2$. Since the polyhedron is regular $\angle C_0C_2C_1 = \pi/p$. Consider the unit sphere around $C_0$ - it is cut by the edges and faces through $C_0$ into spherical polygon with $q$ sides (the vertex figure of the polyhedron). Then $\theta_1$ is half the length of the edge of any such polygon and $\theta_2$ is half the magnitude of its angle. Finally, since $\pi/q$ is half the angle which an edge of the polygon subtends at its center, we can use Napier's rules for a spherical right triangle

\[cos(B)= sin(A)cos(b) \rightarrow cos\frac{\pi}{q} = sin\theta_2 cos\theta_1 = sin\theta_2 sin \frac{\pi}{p}\]

where $2 \theta_2$ is the dihedral angle of the polyhedron. \newline

The results for the dihedral angles of the regular polyhedrons with a given $p$ and $q$ determine how many polyhedrons can be placed around an edge of the constructed polychoron (4-polytope) which is the value of $r$. They are shown in the next table:

\begin{center}
\begin{tabular}{ l l l l }
  p & q & $2\theta_2$ & r \\
  \hline
  3 & 3 & $70.5^o$ & 3, 4, 5 \\
  3 & 4 & $109.5^o$ & 3 \\
  4 & 3 & $90^o$ & 3, 4 \\
  3 & 5 & $138.2^o$ & - \\
  5 & 3 & $116.6^o$ & 3 \\
\end{tabular}
\end{center}

Thus, the only case of polyhedra completely filling a space is the case of 4 cubes around an edge (8 around a vertex) which identifies the case 434 as the only 4-polytope tessellating an Euclidean 3-dimensional space, i.e. the cubic honeycomb which is the only regular honeycomb in 3 dimensions. The regular closed polychora which can be viewed as regular tessellations of a hypersphere can only have 3, 4 or 5 tetrahedra, 3 octahedra, 3 or 4 hexahedra (cubes) or 3 dodecahedra at an edge. They represent the cases 333, 334, 335, 343, 433, and 533 where 333 and 343 are self-dual, and the others are in dual pairs as 433 and 334 and 335 and 533. \newline

Using the two conditions for regularity, the following relations can be obtained

\[N_0qr_q=N_1, \mbox{ } N_1r=N_2p, \mbox{ } N_2 = N_3qq_p\]

where

\[q_p = \frac{2}{2(p+q)-pq}, \mbox{  } r_q=\frac{2}{2(q-r)-qr}.\]

However, they will not be explained further because Euler's equation $N_0 - N_1 + N_2 - N_3 = 0$ is homogeneous in four dimensions and can only determine the ratios between these configurational numbers.

\[N_0 : N_1 : N_2 : N_3 = pq_p : pqq_pr_q : qrp_qr_q : rr_q\]

Therefore, in order to finally determine these polychora, we would have to build them geometrically. Since this gets beyond the scope of this work, here they will simply be identified for the reader:

\begin{itemize}
  \item \{3,3,3\} pentachoron, 5-cell, or 4-simplex (analog of the tetrahedron)
  \item \{4,3,3\} octachoron, 8-cell, or 4-cube (analog of the cube)
  \item \{3,3,4\} hexadecachoron, 16-cell, or 4-orthoplex (analog of the octahedron)
  \item \{3,4,3\} icositetrachoron, 24-cell, or octaplex (no 3-dimensional analog)
  \item \{3,3,5\} hexacosichoron or 600-cell (analog of the icosahedron)
  \item \{5,3,3\} hecatonicosachoron or 120-cell (analog of dodecahedron)
\end{itemize}

One of the algorithms developed for this thesis generates the coordinates of any of these 4-dimensional analogues of the Platonic solids so that various cross-sections can be performed in order to study their symmetry.

\subsection{The regular tessellations in n-1 dimensions and n-polytopes}

As we have seen, the problem of constructing a regular $n$-polytope is only a part of the general problem of constructing a regular tessellation in $n$-1 dimensions. Such a tessellation divides $n$-1-dimensional space into equivalent $n$-1-polytopes which in turn are tessellations of elliptic space in $n$-2 dimensions. In addition, any hypersphere around a vertex of the $n$-1-dimensional tessellation will define a tessellation of elliptic space in $n$-2 dimensions.

\begin{figure}[h!]
  \centering
	\includegraphics[scale=0.3]{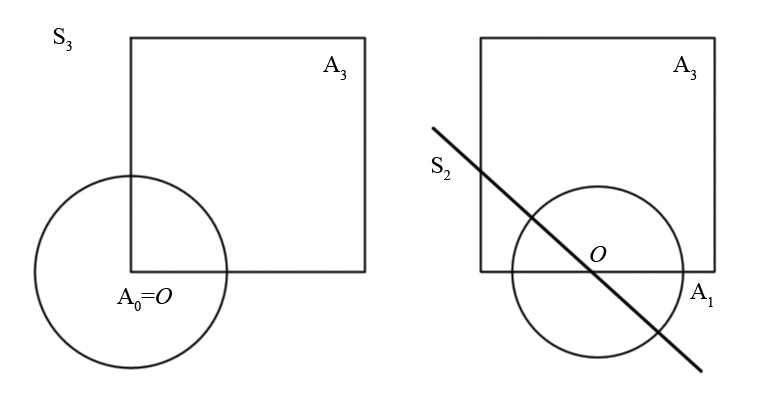}
  \caption{Front view of a cube with cases $p=2$, $q=0$ ($A_q$-vertex, $A_{p+1}$-cube) and $p=2$, $q=1$ ($A_q$-edge, $A_{p+1}-cube$).}
\end{figure}

Consider a $p$+1-facet $A_{p+1}$ and a $q$-facet $A_q$ inside it such that $q \leq p$. In $A_{p+1}$ take any $(p+1-q)$-dimensional flat $S_{p+1-q}$ intersecting $A_q$ at a point $O$ and construct a small hypersphere centered at $O$. It is easy to see that
\begin{flalign*}&\mbox{dim}S_{p+1-q}=p+1-q=\mbox{dim}A_{p+1}-\mbox{dim}A_q\mbox{, }& \\
&\mbox{dim}A_q\cap S_{p+1-q}\leq min(q,p+1-q)\mbox{, and }O\subset A_q\cap S_{p+1-q}& \\
&\Rightarrow_{q+r}A_q\cap S_{p+1-q}=O_r,&\end{flalign*}

i.e. all $q$+$r$-facets through $A_q$ are cut by $S_{p+1-q}$ in $r$-flats through $O$. These flats then cut the hypersphere in $r$-1-dimensional regions (the $q$+1-facets through $A_q$ cut the hypersphere in points, the $q$+2-facets cut it into great circles, etc.) Thus, the figure formed on the hypersphere is a regular tessellation on elliptic space of $p-q$ dimensions. If we denote by $_aN_{bc}$ the number of $b$-facets passing through a $c$-facet and lying on an $a$-facet of the general polytope and by $_aN'_{bc}$ the configurational numbers of the regular tessellation of $p-q$ dimensions, we get

\[_aN'_{bc}=\mbox{}_{a+q}N_{b+q,c+q}\]

When $q=p-1$, a 1-dimensional elliptic regular tessellation can be obtained, i.e. a regular polygon, and $_{p+1}N_{p,p-1}=\mbox{}_2N'_{10}=2$. When $q=p-2$, a 2-dimensional elliptic regular tessellation is obtained as $_{p+1}N_{p,p-2}=\mbox{}_3N'_{20}=\mbox{}_3N'_{10}=\mbox{}_{p+1}N_{p-1,p-2}$. The numbers

\[_{p+1}N_{p,p-2}=\mbox{}_{p+1}N_{p-1,p-2}=k_p \mbox{  } (p = 1, 2, ..., n-1)\]

are analogous to the ones in the 4-dimensional case. For $p=1$, the convention $_rN_{p,-1}=N_{pr}$ is used to denote the total number of $p$-facets in an $r$-facet. \newline

Furthermore, denote the configurational numbers of the $n$-1-facet of the tessellation by $_pF_{qr}$ and the configurational numbers of the vertex figure by $_pV_{qr}$.

\[_{p+1}F_{p,p-2}=\mbox{}_{p+1}N_{p,p-2}=k_p \Rightarrow k_1, k_2, ..., k_{n-2}\]
\[_{p+1}V_{p,p-2}=\mbox{}_{p+2}N_{p+1,p-1}=k_{p+1} \Rightarrow k_2, k_3, ..., k_{n-1}\]

Now, the geometrical method for distinguishing among the elliptic, Euclidean, and hyperbolic case has to be generalized for any dimension. Consider a regular tessellation $k_1k_2k_3...k_n$ of $n$-dimensional space. Let $C_0$ be any vertex, $C_1$ the midpoint of an edge though $C_0$, $C_2$ the center of a face (polygon) through $C_0C_2$ and in general $C_r$ the center of the $r$-facet through $C_0C_1C_2...C_{r-1}$. Then the triangle $C_pC_qC_r$ $(p<q<r)$ is always a right angled triangle with right angle at $Cq$. The angle $C_0C_2C_1$ is half the angle at the center of a plane face subtended by an edge and is equal to $\pi/k_1$. The angle $C_1C_0C_2=\theta_1$ is half the angle between two adjacent edges, $C_2C_1C_3=\theta_2$ is half the dihedral angle between two adjacent plane faces and in general $C_pC_{p-1}C_{p+1}=\theta_p$ is half the dihedral angle between two adjacent $p$-facets. Since there are $k_n$ $n$-1-facets at each $n$-2-facet, the angle $C_{n-1}C_{n-2}C_{n}=\theta_{n-1} = \pi/k_n$. Now consider a sphere at $C_0$ in the 3-facet $C_0C_1C_2C_3$ which is cut by the lines and planes through $C_0$ in a regular spherical polygon with sides $2\theta_1$ and angles $2\theta_2$ and the polygon subtended at the center of the polygon by half the side is $\pi/k_2$. Hence by spherical trigonometry

\[cos\pi/k_2=sin\theta_2 cos\theta_1.\]

Again, in the 4-face $C_0C_2...C_4$ take a hyperplane $H$ perpendicular to $C_0C_1$ at $C_1$ and in this hyperplane consider a small sphere centered at $C_1$. It is cut by the lines and planes in which $H$ cuts the planes and hyperplanes through $C_0C_1$ in a regular spherical polygon with sides $2\theta_2$ and angles $2\theta_3$ and the angle subtended by the center of the polygon by half the side that is $\pi/k_3$. By spherical trigonometry again

\[cos\pi/k_3=sin\theta_3 cos\theta_2.\]

Finally we obtain the formula 

\[cos\pi/k_r=sin\theta_r cos\theta_{r-1}\mbox{ }(r=2,3,...,n-1)\mbox{ with }\theta_{n-1}=\pi/k_n.\]

To determine $\theta_1$, consider the right triangle $C_0C_1C_2$. In Euclidean geometry the sum of its angles $\theta_1+\pi/k_1 + \pi/2 = \pi$ and $\theta_1=\pi/2-\pi/k_1$ while in elliptic geometry it is more and in hyperbolic less. \newline

Now for $n=5$ take regular tessellation in 4-dimensional space. Then

\[cos\pi/k_3=sin\theta_3 cos\theta_2=sin\pi/k_4 cos\theta_2\mbox{ and }cos\pi/k_2=sin\theta_2 cos\theta_1\]
\[sin^2\theta_2+cos^2\theta_2=1 \Rightarrow \frac{cos^2\theta_1 sin^2\theta_2}{cos^2\theta_1}+\frac{sin^2\pi/k_4 cos^2\theta_2}{sin^2\pi/k_4}=1\]
\[\frac{cos^2\pi/k_2}{cos^2\theta_1}+\frac{cos^2\pi/k_3}{sin^2\pi/k_4}=1.\]

If the last expression is greater than 1, the tessellation will be elliptic, and for less than one it will be hyperbolic. \newline

Using these generalization further for $n=5$, we require that the 4-facets and vertex figures of a regular 5-polytope must be elliptic tessellations of 4-dimensional space. Since in four dimensions the only such are the six regular polychora 333, 334, 335, 343, 433, 533 the only possible 5-dimensional cases are

\begin{center}
\begin{tabular}{ l l l }
  3333 & 3334 & 3335 \\
  3343 & 3433 \\
  4333 & 4334 & 4335 \\
  5333 & 5334 & 5335 \\
\end{tabular}
\end{center}

Applying the condition

\[\frac{cos^2\pi/k_2}{cos^2\theta_1}+\frac{cos^2\pi/k_3}{sin^2\pi/k_4}=1\]

we can divide the above cases into three elliptic 3333, 3334, and 4333, three Euclidean 3343, 3433, and 4334 and the other five hyperbolic. Therefore, we can conclude that there are three regular tessellations in an Euclidean space of 4-dimensions, i.e. such space can be completely filled with 8-cells (cubes), 16-cells, or 24-cells and that there are only three regular 5-polytopes - the 5-simplex, the 5-cube, and the 5-orthoplex. \newline

For $n=6$ since the 5-facet $k_1k_2k_3k_4$ and vertex figure $k_2k_3k_4k_5$ must both be elliptic, the tessellations can only be of the form

\[33333, 33334, 43333, 43334\]

which correspond respectively to 6-simplex, 6-cube, 6-othoplex and finally a regular tessellation of a 5-dimensional Euclidean space. It can immediately be seen that in any further dimension the choice will remain the same and the self-dual $n$-simplex $\alpha_n$, the $n$-cube $\gamma_n$ and its dual $n$-orthoplex $\beta_n$ are the only regular polytopes in $n$ dimensions ($n\geq5$) while the $43^{n-2}4$ tessellation $\delta_{n+1}$ is the only regular tessellation of an $n$-dimensional space (notation due to Coxeter \cite{coxeter63}). \newline

\section{The quasi-regular and uniform tessellations and polytopes}

\subsection{Quasi-regular tilings and the crystallographic restriction}

The numbers obtained for each specific regular $n$-polytope in the previous section turn out to be very essential for the description of its properties. They will reappear in their specific combination later in this text as the orders of the generators for Coxeter systems and in the classification of the edges of Coxeter graphs. Here we first introduce the most basic notation which was briefly used above called the Schl\"{a}fli symbol and its extension for quasi-regular tessellations first suggested by Coxeter \cite{coxeter63}. For the sake of brevity this extension is analyzed mainly for the case of polyhedra and respectively tilings and a final extension to uniform tessellations then is introduced with an emphasis on honeycombs and polychora. Nevertheless, quasi-regular honeycombs are briefly introduced in the next section for completeness. Finally, since this subsection will emphasize on tilings, its second part analyses the possible symmetries of all tilings. \newline

The Schl\"{a}fli symbol of a polytope with the numbers $k_1k_2...k_{n-1}$ simply has the form \{$k_1,k_2,...,k_{n-1}$\} and can therefore be used for all derived regular polytopes or tessellations. Therefore, the Schl\"{a}fli symbol of the $n$-1-facet of an $n$-polytope \{$k_1,k_2,...,k_{n-1}$\} is \{$k_1,k_2,...,k_{n-2}$\}, the Schl\"{a}fli symbol of the vertex figure is \{$k_2,...,k_{n-1}$\}, and the Schl\"{a}fli symbol of the dual of the polytope is \{$k_{n-1},k_{n-2},...,k_1$\}. \newline

This notation can now be extended to quasi-regular polyhedra. The interior of the intersection of two dual regular polyhedra \{p,q\} and \{q,p\} centered at the same point has $N_1$ vertices which are exactly the mid-edge points of both \{p,q\} and \{q,p\}. Its faces consist of both $N_0$ \{q\} and \{p\} polygons which are the vertex figures respectively of \{p,q\} and \{q,p\}. There are 4 edges at each vertex and $2N_1$ edges altogether. Then

\[N_0-N_1+N_2=N_1-2N_1+(N_0+N_2)=2\]

and the resulting polyhedron can be denoted as $\left\{\begin{matrix} p\\ q\end{matrix}\right\}=\left\{\begin{matrix} q\\ p\end{matrix}\right\}$. \newline

The possible cases can be derived based on the above restriction:
\begin{flalign*}&p=q=3 \Rightarrow N'_0=N_1=6, N'_1=2N_1=12, N'_2=N_0+N_2=4+4=8& \\
&p=3, q=4 \Rightarrow N'_0=N_1=12, N'_1=2N_1=24, N'_2=N_0+N_2=6+8=14& \\
&p=3, q=5 \Rightarrow N'_0=N_1=30, N'_1=2N_1=60, N'_2=N_0+N_2=12+20=32&\end{flalign*}

For the three elliptic cases we then get
\[\left\{\begin{matrix} 3\\ 3\end{matrix}\right\}=\{3,4\} - \mbox{octahedron}, \left\{\begin{matrix} 3\\ 4\end{matrix}\right\} - \mbox{cuboctahedron, and} \left\{\begin{matrix} 3\\ 5\end{matrix}\right\} - \mbox{icosidodecahedron}\]

and since the edges are all alike, each separating $p$ from $q$, this gives rise to the definition of a quasi-regular polyhedron as a polyhedron with exactly two kinds of regular faces that is edge-transitive. The proof that these are the only quasi-regular polyhedra originates from the fact that the dihedral angles at a vertex make a total that must conform to the inequality

\[r(1-\frac{2}{p})\pi+r(1-\frac{2}{q})\pi < 2\pi.\]

Therefore,
\[1-\frac{2}{p}+1-\frac{2}{q} < \frac{2}{r} \Rightarrow 1-\frac{1}{p}-\frac{1}{q}<\frac{1}{r} \Rightarrow \frac{1}{p}+\frac{1}{q}+\frac{1}{r}>1\]

and since $p$ and $q$ cannot be less than 3, $r=2$ and $p=3,q=4$ or $p=3,q=5$. In addition we can consider the dual pair of regular tilings \{3,6\} and \{6,3\} where the trihexagonal tiling $\left\{\begin{matrix} 3\\ 6\end{matrix}\right\}$ is obtained with vertices as the intersections of the their edges. \newline

The symmetry group of such a tiling is an infinite group of congruent transformations in the plane. This group contains a finite subgroup of index 2. A classification of the symmetries of the plane tilings then leads to the crystallographic restriction theorem which states that if a discrete group of rotations in the plane has more than one center of rotation, then the only rotations that can occur are of order 1, 2, 3, 4, and 6. \newline

Consider the Euclidean motion group $\mathbb{R}^2.O(2)$ of isometries on the plane. Any finite subgroup of this group fixes a point and so is conjugate to a finite subgroup of $O(2)$ that fixes the origin. The finite subgroups of $O(2)$ are then the cyclic groups of order $n$, i.e. rotations of $2\pi/n$, and the dihedral groups of order 2$n$ with rotations as a subgroup of index 2 and reflections that conjugate those to their inverse rotations. The subgroup $\mathbb{R}^2$ consists of all translations. Consider all discrete cases, i.e. rotations and translations that cannot be arbitrarily close to the identity transformation and are bounded from below. Let $\Gamma \leq O(2)$ be discrete subgroup and consider the lattice group $L=\Gamma \cap \mathbb{R}^2=\mathbb{Z}a+\mathbb{Z}b$ generated from translations by two linearly independent vectors. Furthermore, let $\bar{\Gamma}$ be the image of $\mathbb{R}^2.O(2)$ in $\mathbb{R}^2.O(2)/\mathbb{R}^2$, i.e. $\bar{\Gamma}=\Gamma/L$. \newline

Then $\bar{\Gamma}$ preserves the lattice $L$. The proof is the following. Take a vector $b \in L$ or equivalently a translation by this vector $t_b \in \Gamma$ and take $\gamma \in \Gamma$ that maps to $\bar{\gamma} \in \bar{\Gamma}$. In terms of linear operators set
\begin{flalign*}&\gamma(v)=Av, \mbox{ } t_b(v)=v+b, \mbox{ } \gamma^{-1}(v)=A^{-1}v& \\
&\gamma t_b \gamma^{-1}(v)=\gamma t_b(A^{-1}v)=\gamma(A^{-1}v+b)=A(A^{-1}v+b)=v+A(b)&\end{flalign*}

Therefore $\gamma t_b \gamma^{-1}=t_{\bar{\gamma}(b)}$ is a conjugate translation by $\bar{\gamma}(b)$ and since $\gamma, t_b \in \Gamma$ by the operation closure it follows that $t_{\bar{\gamma}(b)} \in \Gamma$ and finally $\bar{\gamma}(b) \in L$. \newline

Because of this established fact, $\bar{\Gamma}=C_n$ or $\bar{\Gamma}=D_{2n}$ with $n=1,2,3,4,6$ and $\bar{\Gamma}=C_n$ as rotation parts has maximum order of 12. A proof in terms of linear operators is as follows. Let $A \in \bar{\Gamma}$ be a rotation, i.e. $detA=1$. We have to show that the order of $A$ is 1,2,3,4 or 6. Consider the characteristic polynomial of $A$ $x^2-tr(A)x+det(A)=x^2-tx+1$. Consider $A$ to be a rotation by $\theta$. Then since rotation is analogical to scaling with complex numbers, $x^2-tx+1$ has complex roots, i.e. $t^2-4\leq0$ (the only real cases are $+1$ and $-1$). Thus, the matrix is diagonalizable over the complex numbers and $tr(A)=t=z+\bar{z}$ is a real number. Furthermore, since $A$ stabilizes the lattice $L=\mathbb{Z}a+\mathbb{Z}b$, it takes both $a$ and $b$ to integer multiples of $a,b$. If we take these as column vectors of $A$ in the basis $a,b$, $A$ will have integer entries in this basis and thus its trace $t$ is an integer. Then $t^2-4 \leq 0 \Rightarrow t = \pm1,\pm2, 0$ where $t=\pm2=2cos(\theta)$ are the rotations of order 1 and 2, $t=\pm1=2cos(\theta)$ are the rotations of order 3 and 6, and $t=0=2cos(\theta)$ are the rotations of order 4. This concludes the proof. \newline

Considering Gaussian integers for coordinates, the symmetry group of the tiling \{4,4\} is then generated by the translation $z'=z+1$ and the rotation $z'=iz$ of order 4 while the symmetry group of the \{3,6\} is generated by the same translation along with a rotation $z'=e^{\pi i/3}z=(e^{2\pi i /3}+1)z$ of order 3. The crystallographic restriction theorem can be used to classify all symmetries of a tiling on the plane with translations, rotations, reflections, and glide reflections as isometries of the Euclidean plane, also called wallpaper groups or crystallographic groups on the plane.

\subsection{Uniform honeycombs and geometric operations}

The extension of the Schl\"{a}fli symbol introduced by Coxeter can be generalized for rectified $n$-polytopes where rectification is the process of taking the intersection of two dual polytopes, i.e. cutting the vertices at the midpoints of the edges which will result in a polytope bounded by both the vertex figures and rectified faces of the original. \newline

For a quasi-regular honeycomb, all cells must be regular and all vertex figures must be quasi-regular. Alternative conditions then are that the vertex figures are all the same and the cells are of two alternating kinds. The only two regular polyhedra, whose angle sum divides $2\pi$ are the octahedron and the tetrahedron (their sum is $\pi$). The only quasi-regular honeycomb then is $\left\{\begin{matrix} 3, 3\\ 4\end{matrix}\right\}$ or the alternated cubic honeycomb. It can be seen as a cubic honeycomb with alternate vertices removed reducing cubic cells to tetrahedra and creating octahedron cells in the remaining gaps. \newline

Rectification is a special case of the more general truncation operation which can be used to derive a list of uniform polytopes and tessellations. Although there is a new notation known as the Wythoff symbol originating from the Wythoffian construction of uniform polytopes, there is also a further extension of the Schl\"{a}fli symbol and a notation that will be introduced in the next section that adds more information about the symmetries of the tessellation (resp. of the polytope). The final extended Schl\"{a}fli symbol denotes the $k$-th rectification of a polytope as $t_k\{p_1, p_2, ..., p_{n-1}\}$ where $t_{0,1}$ is a truncation applied to polygons or polytopes of higher dimension, $t_{0,2}$ is cantellation (both edges and vertices removed) applied to polyhedrons or higher, $t_{0,3}$ is runcination, $t_{0,1,2}$ is cantitruncation (cantellation and trucation), $t_{0,1,2,3}$ is runcicantitruncation, etc. \newline

All these geometric operations can be generalized in terms of sequences of cross-sections with higher dimensional space. Truncation results from the sequence of cross-sections parallel to a facet of the vertex figure of the polytope where crossing the vertex leads to intersection of the previously crossed edges and origination of the next. Cantellation is analogously the result of sequence of cross-sections parallel to an edge of the polytope where crossing the edge leads to intersection of previously crossed faces and origination of the next. Such a generalization can be used for the reconstruction of an $n+1$-dimensional polytope through the geometrical operations performed on its cross-sections.\newline  

Using this, we were able to reconstruct the graph of a 4-polytope with cross-sections which are the gradual truncations of a cube to its dual. \newline

These operations (besides rectification) introduce the semi-regular polytopes in addition to the quasi-regular and regular ones. All of them are unified under the definition of a uniform polytope whose only two conditions are: (1) uniform polytope facets; (2) vertices of the same kind (vertex-transitivity). The uniform 2-polytopes are necessarily the regular polygons. \newline

There are 11 uniform tilings of the plane and 28 uniform convex honeycombs in 3 dimensions, also called the Archimedean honeycombs. From the latter, there are just one regular (cubes) and one quasi-regular (octahedra and tetrahedra), both mentioned above. Truncation has been used to derive 7 additional honeycombs originating from the cubic one and 4 additional originating from the alternated cubic honeycomb. Finally there are 15 more from prismatic forms derived from modifications of the uniform plane tilings. \newline

\section{Coxeter groups and reflection groups}

\subsection{Coxeter systems}

In the previous section it was shown that the group of isometries on the Euclidean plane admits only a few discrete subgroups of $O(2)$ thus allowing a classification of the possible symmetries of a tiling that is not necessarily uniform. The 3-dimensional case of crystallographic groups is the space groups. A formal description of the symmetries of a tessellation (resp. polytope) is the goal of this section which leaves the techniques of $n$-dimensional geometry aside and uses tools mainly from abstract algebra. \newline

A Coxeter system is a pair $(W,S)$ where $W$ is a group and $S$ a set of generators $S\subset W$ restricted by the relations $(s_is_j)^{m(s_i,s_j)}=1$ where $m(s_i,s_i)=1$, i.e. $s_i$ is an involution, and $m(s_i,s_j)=m(s_j,s_i) \geq 2, \mbox{ } \forall s_i \neq s_j \in S$. If there is no relation between $s_i$ and $s_j$ the convention $m(s_i,s_j)=\infty$ is used. The group $W$ is then the quotient $F/N$ where $F$ is a free group on the set $S$ and $N$ is the normal subgroup generated by all elements $(s_is_j)^{m(s_i,s_j)}$. Furthermore, $|S|=n$ is the rank of the Coxeter system. $W(M)$ is then a Coxeter group that can be constructed from a symmetric $n \times n$ matrix $M = (m_{ij})_{1\leq i,j\leq n}$ indexed by $S$ with entries in $\mathbb{Z}\cup \{ \infty \}$ such that $m_{ii}=1$ and $m_{ij}\geq2 \mbox{ } \forall i\neq j$. The Coxeter group of type $M$ is then analogously

\[W(M)=\langle s_1,s_2,...\in S| (s_is_j)^{m_{ij}}=1,m_{ij}\in M\rangle\]

which will be denoted $W$ when no ambiguity is possible. \newline

A few lowest rank examples are the following. If $|S|=n=1$ then $M=(1)$ and $W(M)=\langle s_1|s_1^2=1\rangle$ which is the cyclic group of order 2. For $n=2$ then $M=\bigl(\begin{smallmatrix}1&m\\m&1\end{smallmatrix} \bigr) \Rightarrow W(M)=\langle s_1,s_2|s_1^2=1,s_2^2=1,(s_1s_2)^m=1\rangle=D_{2m}$ for $m\in\mathbb{N}\cup\infty$ which is the Klein Four group for $m=2$, the dihedral group of finite order for $2\leq m\leq \infty$ and the infinite dihedral group $D_{\infty}$ for $m=\infty$. \newline

Instead of the matrix, the Coxeter system (W,S) can be constructed from an undirected graph $\Gamma$ with a vertex set $S$ where two vertices $s$ and $s'$ are joined with an edge that is labeled $m(s,s')$ if $3 \leq m \leq \infty$. Therefore if the distinct vertices $s$ and $s'$ are not joined, then $m(s,s')=2$. The edges with label $m(s,s')=3$ are omitted due to their frequency and by convention. The resulting graph $\Gamma$ is called Coxeter graph and as a notation contains more information than the extended Schl\"{a}fli symbol. \newline

\begin{figure}[h!]
  \centering
    \includegraphics[width=0.4\textwidth]{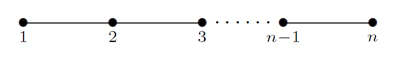}
  \caption{Sample Coxeter graph.}
\end{figure}

The entire information about the matrix can be reconstructed from the Coxeter graph as the same applies to the Schl\"{a}fli symbol. To show this we need to study the connection of the Coxeter groups to the reflection groups and the groups of symmetries of the regular polytopes. \newline

Because Coxeter groups are based on an abstract presentation, they do not necessarily admit a faithful representation as reflection groups. The abstract group of each reflection group is a Coxeter group as the reflections are a special case of involutions \cite{coxeter63}. Infinite Coxeter groups in particular may not admit a representation as a reflection group. However, finite Coxeter groups have a faithful linear representation as groups generated by reflections in Euclidean space. Some of the finite linear groups generated by such reflections are groups of symmetries of the regular polytopes in Euclidean space. We will show that Coxeter groups admit a representation since this is essential for the purpose of this work. \newline

Let $(W,S)$ be a Coxeter system of type $M$ and take $n=|S|$. We will now construct a real linear representation of $W$ of degree $n$ such that the images of the elements of $S$ are reflections in $\mathbb{R}^n$. A reflection on a real vector space $V$ is a linear transformation on $V$ fixing a subspace of $V$ of codimension 1, a reflection hyperplane $H_\alpha$, and having a nontrivial eigenvector $\alpha$ with eigenvalue $\lambda=-1$, called a root of the reflection. Now take the vector space $V$ over the field $\mathbb{R}$ with a basis $\{\alpha_s|s\in S\}$. Define a bilinear form $B$ on $V$ such that

\[B(\alpha_s,\alpha_{s'})=-cos\frac{\pi}{m(s,s')}\]

where $m(s,s')=m_{ij}\in M$ and since $m_{ij}=m_{ji}$ the form is symmetric. Then $B=-1$ for $m(s,s')=\infty$, $B(\alpha_s.\alpha_s)=1$, and $B(\alpha_s.\alpha_s)\leq 0$ for $s\neq s'$ with equality only for $m(s,s')=2$. Now for each $s \in S$ define a reflection $\sigma_s:V\rightarrow V$ such that $\sigma_s\lambda=\lambda-2B(\alpha_s,\lambda)\alpha_s$. Then $\sigma_s\alpha_s=-\alpha_s$ and $\sigma_s\lambda=\lambda$ with $\{\lambda \in V|B(x,\lambda)=0\} := H_s$ the hyperplane orthogonal to $\alpha_s$. Therefore $\sigma_s$ has order 2 in $GL(V)$. Additional observation from the definition is that $B(\sigma_s\lambda,\sigma_s\mu)=B(\lambda,\mu) \mbox{ } \forall \lambda, \mu \in V$, i.e. the reflection preserves the bilinear form and each element generated by $\sigma_s$ will preserve it. A final observation is that $|\sigma_s\sigma_{s'}|=m(s,s')$ and should be proven next. \newline

Consider $V_{s,s'}:=\mathbb{R}\alpha_s\oplus\mathbb{R}\alpha_{s'}$. The restriction of $B$ to $V_{s,s'}$ is positive semidefinite since for $\lambda=a\alpha_s+b\alpha_{s'}$ we get
\begin{flalign*}&B(\lambda,\lambda)=B(a\alpha_s+b\alpha_{s'},a\alpha_s+b\alpha_{s'})=& \\
&=a^2B(\alpha_s,\alpha_s)+2abB(\alpha_s,\alpha_{s'})+b^2B(\alpha_{s'},\alpha_{s'})=& \\
&=a^2-2abcos(\pi/m(s,s'))+b^2=& \\
&=a^2-2abcos(\pi/m(s,s'))+b^2(cos^2(\pi/m(s,s'))+sin^2(\pi/m(s,s'))=& \\
&=(a-bcos(\pi/m(s,s')))^2+b^2sin^2(\pi/m(s,s'))\geq0.&\end{flalign*}

Therefore the form is positive definite if $sin(\pi/m(s,s'))\neq0$ and $m<\infty$. Note further that $\sigma_s$ and $\sigma_{s'}$ stabilize $V_{s,s'}$, so the order of $\sigma_s\sigma_{s'}$ as an operator in $V_{s,s'}$ has two cases. (1) If $m<\infty$ since the form is positive definite we can consider Euclidean plane and since $B(\alpha_s,\alpha_{s'})=-cos(\pi/m(s,s'))=cos(\pi-(\pi/m(s,s')))$ and the angle between $\mathbb{R}\alpha_s$ and $\mathbb{R}\alpha_{s'}$ is therefore $\pi-(\pi/m(s,s'))$, the angle between $H_\alpha$ and $H_{\alpha'}$ is $\pi/m(s,s')$. Because rotation through $2\pi/m(s,s')$ can be achieved as a product between two reflections with an angle between their fixed hyperplanes $\pi/m(s,s')$, it follows that $\sigma_s\sigma_{s'}$ has order $m(s,s')$. The fact that $B$ is positive definite on $V_{s,s'}$ implies that $V=V_{s,s'}\oplus V_{s,s'}^{\perp}$ where $V_{s,s'}^{\perp}:=s^{\perp} \cap s'^{\perp}$ is fixed by both $\sigma_s$ and $\sigma_{s'}$. Then $\sigma_s\sigma_{s'}$ has order $m(s,s')$ also on $V$. (2) If $m=\infty$, $B(\alpha_s,\alpha_{s'})=-1$. Then if $\lambda=\alpha_{s}+\alpha_{s'} \mbox{, } B(\lambda,\alpha_s)=B(\lambda,\alpha_{s'})=0$ and $\sigma_s$ and $\sigma_{s'}$ fix $\lambda$. Then $\sigma_s\sigma_{s'}\alpha_s=\sigma_s(\alpha_s+2\alpha_{s'})=3\alpha_s+2\alpha_{s'}=2\lambda+\alpha_s$. Applying $\sigma_s\sigma_{s'}$ $k\in\mathbb{Z}$ times, $(\sigma_s\sigma_{s'})^k\alpha_s=2k\lambda+\alpha_s$. Therefore we can conclude that $\sigma_s\sigma_{s'}$ has infinite order on $V_{s,s'}$ and therefore also on V.\newline

This helps us conclude that there is a unique homomorphism $\sigma:W\rightarrow GL(V)$ where $\sigma(s_1s_2...s_n)=\sigma(w)=\sigma_w=\sigma_{s_1}\sigma_{s_2}...\sigma_{s_n}$. We call this homomorphism a linear representation of $W$. The fact that this representation is then faithful is a corollary from another theorem regarding the length function for the Coxeter group and is outside of the scope of this text. However, one final definition which is necessary for the later sections is as follows. The subgroup $W_I$ generated by a given subset $I\subset S$ and any of its conjugates are called parabolic subgroups of the Coxeter group $W$ and are in fact Coxeter groups themselves.

\subsection{Root systems}

Consider a Coxeter system $(W,S)$ of type $M$. The set $\Phi :=\bigcup_{s \in S, w \in W} \sigma(w)\alpha_s$ is called the root system of $W$. It consists of the collection of orbits of the unit vectors in the reflection representation space $V$ on which $W$ acts, i.e. the collection of unit vectors in $V$ permuted by $W$. These remain unit vectors because $W$ preserves the bilinear form. Furthermore, $\Phi=\Phi^+\cup\Phi^-$ where

\[\Phi^+:=\Phi\cap(\sum_{s \in S}\mathbb{R}_{\geq0}\alpha_s)\mbox{, }\Phi^-:=\Phi\cap(\sum_{s \in S}\mathbb{R}_{\leq0}\alpha_s).\]

A root $\alpha$ is called positive (write $\alpha>0$) if $\alpha\in\Phi^+$ and negative ($\alpha<0$) if $\alpha\in\Phi^-$. Note that $\Phi=-\Phi$ since $\sigma_s(\alpha_s)=-\alpha_s$, $\Phi\cap\mathbb{R}\alpha=\{\alpha,-\alpha\} \mbox{ } \forall\alpha\in\Phi$, and $\sigma(W)\Phi=\Phi$. Positive and negative roots are defined with regard to certain total ordering like lexicographic ordering where $\sum a_i\lambda_i<\sum b_i\lambda_i \Rightarrow a_k < b_k$ where $k$ is the least index such that $a_k\neq b_k$.\newline

A subset $\Delta \subset \Phi^+$ of vectors $\alpha_i$ constituting a basis for the $\mathbb{R}$-span of $\Phi$ in $V$ is called a simple system. A root $\alpha \in \Delta$ is called simple root of $\Phi$. Each $\alpha' \in \Phi$ is a linear combination of $\Delta$ with coefficients all of the same sign. \newline

Further relation between the roots and the reflections can be established through the consideration of the set $R=\{wsw^{-1}|w\in W,s \in S\}$ or reflections of the Coxeter system $(W,S)$. By the geometric representation $\sigma : W \rightarrow GL(V)$, each $s\in S$ acts on $V$ as a reflection $\sigma(s)\alpha_s$. More generally, a reflection in $GL(V)$ can be associated to each root $\alpha \in \Phi$. Consider $\alpha=\sigma(w)\alpha_s := w(\alpha_s)$ (for brevity) for $w\in W \mbox{ and } s\in S$. Then $wsw^{-1}\in R$ acts on $W$ as follows:
\begin{flalign*}&\sigma(wsw^{-1})\lambda=:wsw^{-1}(\lambda)=w[w^{-1}(\lambda)-2B(w^{-1}(\lambda),\alpha_s)\alpha_s]=& \\
&=\lambda-2B(w^{-1}(\lambda),\alpha_s)w(\alpha_s)=\lambda-2B(\lambda,w(\alpha_s))w(\alpha_s)=& \\
&=\lambda-2B(\lambda,\alpha)\alpha&\end{flalign*}

which shows that $wsw^{-1}$ does not depend on the choice of $w$ and $s$ but only on the choice of $\alpha$ so can be denoted $s_\alpha$. Furthermore, $s_\alpha$ acts on $V$ as a reflection sending $\alpha$ to $-\alpha$, fixing pointwise the hyperplane $H_\alpha \perp \alpha$. As a consequence, $\alpha$ and $-\alpha$ both determine the same reflection $s_\alpha=-s_\alpha$. The root-reflection correspondence is finally established due to the bijective map $\alpha \rightarrow s_\alpha$ (for $\alpha\in\Phi^{+}$). Therefore, each reflection of a Coxeter group $W$ has a unique positive root $\alpha\in\Phi^+$ and each $\alpha\in\Phi^+$ is the positive root of a unique orthogonal reflection with respect to $B$. If $\alpha,\beta\in\Phi$ and $w,w'\in W$ such that $w'(\beta)=w(\alpha)$, then $ws_\alpha w^{-1}=w's_\beta w'^{-1}$. \newline

This correspondence helps us interpret the relations in a Coxeter system of the form $(s_i s_j)^{m_{ij}}$ as the result of two reflections fixing hyperplanes meeting at an angle $\pi / m_{ij}$. The element $s_i s_j \in S$ being of order $m_{ij}$ than has the geometrical interpretation of a rotation by $2\pi/m_{ij}$. \newline

Some further important properties and definitions of root systems to be mentioned here are as follows. If $\Delta$ is a simple system in $\Phi$, then $(\alpha,\beta)\leq0\mbox{ }\forall\alpha\neq\beta\in\Delta$. If $\alpha_s\in\Delta$ then $s(\Pi\backslash\{\alpha_s\})=\Pi\backslash\{\alpha_s\}$. A root system $\Phi$ is crystallographic if it satisfies the additional requirement

\[\frac{2(\alpha,\beta)}{(\beta,\beta)}\in\mathbb{Z}\mbox{ }\forall\alpha,\beta\in\Phi.\]

The vector $\alpha^\vee:=2\alpha/(\alpha,\alpha)$ is called a coroot of $\alpha\in\Phi$. The set $\Phi^\vee$ of all coroots is the dual root system or the inverse root system of $\Phi$. The reflections $w\in W$ generated by $\Phi^\vee$ are the same as $\Phi$, i.e. $w\alpha^\vee=w(\alpha)^\vee$. \newline

A final definition for this section should be given for the fundamental domain of the action of $W$ on $V$. Take a positive system $\Phi^+$ containing a simple system $\Delta$ and consider the open half-spaces $H_\alpha^+:=\{\lambda\in V|(\lambda,\alpha)>0\}$ and $H_\alpha^-:=-H_\alpha^+$ of each reflecting hyperplane $H_\alpha$. Define $C:=\cap_{\alpha_s\in\Delta}H_\alpha^+$ which is open and convex as intersection of open and convex sets. Let $D:=\bar{C}$ be the closure of $C$, i.e. the intersection of the closed half-spaces $H_\alpha^+\cup H_\alpha$. Then

\[D = \{\lambda\in V|(\lambda,\alpha)\geq 0 \mbox{ } \forall \alpha\in\Delta\}\]

and each $\lambda\in V$ is conjugate to exactly one point in $D$. Thus $D$ is called a fundamental domain for the action of $W$ in $V$. \newline

This could offer a different insight into the geometric representation of the Coxeter group. The nodes of the Coxeter graph represent the walls of the fundamental domain and two nodes are joined by a branch whenever the corresponding walls are not perpendicular. Moreover, the branches are marked with numbers $m_{ij}>2$ to indicate the angles $\pi/m_{ij}$. In the case of a connected graph without any even marked branches, all the reflections in the group are conjugate to one another \cite{coxeter63}. This interpretation of the Coxeter graphs in terms of fundamental domains originally proposed by Coxeter leads to their final extension to the geometrical operations discussed in the previous section and the uniform tessellations and polytopes.

\subsection{Classification}

All symmetry groups of regular polytopes are finite Coxeter groups (and resp. finite reflection groups). All symmetry groups of regular tessellations are affine Coxeter groups (and resp. affine reflection groups containing normal abelian subgroups such that the quotient group is finite and is itself a Coxeter group). The Coxeter graph of an affine Coxeter group is obtained by adding an additional vertex as in the construction performed in section 2. The dual polytopes or tessellations have the same symmetry groups and therefore the same Coxeter groups. \newline

For the regular polytopes in any number of dimensions, the symmetry groups are respectively the symmetric group $S_{n+1}$ or the Coxeter group of type $A_n$ for the regular $n$-simplex $\alpha_n$, and the hyperoctahedral group or the Coxeter group of type $B_n=C_n$ for the $n$-cube $\gamma_n$ and its dual the $n$-orthoplex $\beta_n$. The root systems for these use the same notation with the difference that $B_n$ and $C_n$ are interchanged (dual root systems with $B_2$ and $C_2$ isomorphic). The root system $B_n$ of $\gamma_n$ which is specifically important for one of the algorithms described in the next section is constructed as follows. \newline

Let $V = R^n$, and let $\Phi$ consist of all integer vectors in $V$ of length 1 (short roots) or $\sqrt{2}$ (long roots). The total number of roots is $2n^2$ with $2n$ short roots $\pm e_i$ and $2n(n-1)$ long roots $\pm e_i\pm e_j (i<j)$. For $\Delta$ take the simple long roots $\alpha_i = e_i - e_{i+1}$, for $1\leq i\leq n-1$, and the short root $\alpha_n=e_n$. The reflection $s$ through the hyperplane perpendicular to the short root $\alpha_n$ is then the negation of the $n$-th coordinate. \newline

Furthermore, the symmetry group of the pentagon is $H_2$, the symmetry group of the dodecahedron and its dual icosahedron is the full icosahedral group $H_3$ and the symmetry group of their 4-dimensional analogues (the 120-cell and the 600-cell) is $H_4$. The symmetry group of the 24-cell is $F_4$. The Coxeter groups of type $D_n$ ($n$-demihypercube), $E_6$ ($2_{21},1_{22}$), $E_7$ ($3_{21}, 2_{31}, 1_{32}$), and $E_8$ ($4_{21}, 2_{41}, 1_{42}$) are the symmetry groups of certain semiregular polytopes. The symmetry group of the hexagon is $G_2$.\newline

The affine Coxeter groups are then classified as $\tilde{A_n}$ for the simplectic uniform tessellation, $\tilde{B_n}$ for the demihypercubic uniform tessellation, $\tilde{C_n}$ for the hypercubic uniform tessellation, $\tilde{E_6}$ for $2_{22}$, $\tilde{E_7}$ for $3_{31},1_{33}$, $\tilde{E_8}$ for $5_{21},2_{51},1_{52}$, $\tilde{F_4}$ for 16-cell and 24-cell uniform tessellations and $\tilde{G_2}$ for hexagonal and triangular tiling. \newline

The final notation to be reviewed in this work is the ringed Coxeter graph which contains enough information to explicitly enumerate almost all types of uniform polytopes and uniform tessellations. Each uniform polytope can be generated using the mirror hyperplanes bounding the fundamental region and a single generator point. The reflections of the point through the mirror hyperplanes and their further reflections through the same hyperplanes form the set of vertices of the polytope. The edges of the polytope connect each point to a mirror point; the faces can be constructed as cycles of edges, etc. The location of the generating vertex is specified as all nodes of the Coxeter graph corresponding to the mirror hyperplanes on which the vertex does not lie are marked with a ring (equidistant from all ringed node hyperplanes). Thus, all mirror hyperplanes where the generating vertex lies do not generate additional vertices. A diagram needs at least one active node to represent a polytope and therefore all Coxeter graphs of the regular polytopes have their first node ringed. \newline

The more general case of uniform polytopes with one marked hyperplane corresponds to a generating point at a vertex of the fundamental domain (which is always a simplex because of the way it is constructed). If all nodes are ringed, the generator point lies in the interior of the simplex. Generally, if $n$ nodes of the Coxeter graph are marked, the generating point gets $n-1$ degrees of freedom on $n$-1-facets of the fundamental domain and can generally be taken to be at the center the $n$-1-facet for equal $n$-1-faces of the final polytope. A secondary feature can be used for the special cases of uniform polytopes with non-reflectional symmetry where the central dot of a ringed node is removed to imply alternate nodes deleted. The constructed polytope will then only have a subsymmetry of the original Coxeter group. Eventually if all nodes are marked in this way, the polytope is called a snub. Using this final notation we can describe for example the cuboctahedron, rhombicuboctahedron, octahedron, truncated cube and other uniform polytopes derived from geometric operations on the cube.

\section{The hypercube tessellation and convex polytope algorithms}

\subsection{The developed algorithms}

Two of the algorithms developed for this work will be presented in this section. Since they consist of some common approaches and differ by removable components, they will be described in terms of one sequentially consistent algorithm. Although our initial plan was to use parametric equations, it offered no good solution when it comes to the calculation of the intersection. Furthermore, algebraic equations did not allow for plotting of vertical lines. \newline

The hypercube tessellation algorithm generates an $n$-cube and $n$ additional $n$-cubes on all of its sides and obtains a cross-section with a 3-dimensional space parallel to three independent vectors (preferably roots from the $B_n$ root system). This results in the 3-dimensional plot of a few space-filling polyhedra that build up the tessellation and an additional optional 2-dimensional plot. The convex polytope algorithm generates a convex $n$-polytope by analyzing its possible $k$-facets from its vertex coordinates and then performs a similar cross-section through specified parameters. \newline

For an $n$-dimensional convex polytope with $V$ as the set of its vertices and $F$ as the set of its $n$-facets, take all subsets $S \subset V, |S|=n$. For these $n$ points $A_1, A_2,...,A_n$ find a nonzero vector $(u_1,u_2,...,u_n)$ such that it is orthogonal to the vectors $A_i-A_n, i = 1, 2, ..., n-1$. Thus, we need to solve a homogeneous system of $n-1$ equations in $n$ variables $u_1,u_2,...,u_n$. The nullspace of the constructed matrix can be more than one dimension in case the $n$ vertices constitute a subspace of the hyperplane (there at least two linearly independent non-zero vectors in the kernel of the considered linear map). These cases immediately imply that the given $S$ is not an $n$-1-facet of the polytope and therefore can be ignored. \newline

The remaining results are the normal vectors $\vec{n_j}$ of the hyperplanes $H_j$ passing through the points from $S$ where $j$ indexes each of the remaining cases for $S$. It is important to consider both positive and negative orientation of the $\vec{n_j}$ when $H_j$ is tested for an $n$-face. For this purpose one of the definitions of convexity for the polytope has been used. \newline

A polytope is convex if it lies entirely on one side of each of its $n$-1-facets. Therefore, if the polytope is convex and has a facet lying in the hyperplane $H_j$, all points $A \in V\backslash S$ must lie in the closed negative half-space $H_j^-$. The signed distance is obtained from the projection of any vector from the hyperplane to the point $A=(a_1,a_2,...,a_n)$ onto the normal vector of the plane $\vec{n_j}=(u_1,u_2,...u_n)$.
\begin{flalign*}&\frac{|(\vec{a}-\vec{x})\cdot \vec{n_j}|}{|\vec{n_j}|}=\frac{|a_1u_1+a_2u_2+...+a_nu_n-u_1x_1-u_2x_2-...-u_nx_n|}{|\vec{n_j}|}=& \\
&=\frac{|a_1u_1+a_2u_2+...+a_nu_n-(-C)|}{|\vec{n_j}|} = \frac{a_1u_1+a_2u_2+...+a_nu_n+C}{\sqrt{u_1^2+u_2^2+...+u_n^2}}&\end{flalign*}

The translation constant $C$ can be obtained from the regular hyperplane equation of $H_j$ evaluated at any point $A \in S$. The convex hull defined by the vertices in $V$ is then

\[ \bigcap H_{j_k}^- = conv(V), k = 0, 1,...,|F|\]

which is the desired polytope and the solution of $|F|$ inequalities. \newline

However, only the normal vectors $n_{j_k}$ of the $n$-facets of the polytope are necessary for the second part of the algorithm which is the cross-section with a 3-subspace where we can observe certain symmetries of the polytope. Three linearly independent vectors and a translation point are sufficient for defining a unique 3-subspace inside the n-space. In order to observe symmetries of the polytope in the 3-subspace, only specific orientations are allowed. These are determined by the root system $\Phi$ of the Coxeter group for the specific polytope. Thus, picking a 3-subspace which is parallel to three linearly independent roots (not necessarily simple) ensures that the observed cross-section conforms to the symmetries of a parabolic subgroup of the original Coxeter group. It then can in term be described with the root system stabilized by this subgroup. For example considering the root system $B_4$ and a hypercube, we can pick the long roots

\begin{center}
\begin{tabular}{ r l }
  $u_1$ = (-1,-1,0,0) & $\angle(u_1,u_2)=\pi/3$ \\
  $u_2$ = (-1,0,-1,0) & $\angle(u_2,u_3)=\pi/4$ \\
  $u_3$ = (-1,0,0,0) & $\angle(u_3,u_1)=\pi/4$ \\
\end{tabular}
\end{center}

which will result in $B_3$ root system and the cross-section with the hypercube should be invariant under all reflections along these roots. The example cross-section is a cube standing on its vertex with respect to the xy-plane and the resulting 2-dimensional cross-section is a hexagon. Therefore, instead of picking any possible orientation of the 3-subspace which is also possible, we pick only vectors $\alpha \in \Phi$. \newline

It is important to note that the choice of a normal vector of the 3-subspace is not unique if its codimension in the $n$-space is greater than 1. For a larger codimension, i.e. larger dimension of the orthogonal complement of the 3-subspace, we would have to select a basis of $n-3$ vectors in order to uniquely identify the 3-subspace. Therefore, a much better approach is to simply select three vectors spanning the 3-subspace in order to determine its orientation in the $n$-space. \newline

After the three roots are selected, the Gram-Schmidt orthogonalization process can be used to produce the new orthonormal basis and a transition matrix $T$ with the unit vectors as columns vectors. Once the linearly independent roots are chosen, orthogonalization inside the 3-subspace naturally does not influence the cross-section in any way. The $n$-facet normal vectors $T^{-1}n_{j_k} = T^tn_{j_k}$ then define new positive half-spaces which reorients the convex polytope. In the case of the hypercube tessellation algorithm, the same transformation is applied to all side-cubes with the difference that they are initially translated along some initial unit vector. When the half-space inequalities are produced, translation of the hyperplanes is applied using the point specifying the translation of the 3-subspace. Finally, the inequalities are solved for $x_1, x_2,$ and $x_3$ given $x_i = 0, 3 < i \leq n$ and the result is ready for plotting. \newline

\subsection{Obtained results}

The following subsection introduces some important observations and results obtained in this work from both the algorithms and the performed study. \newline

Consider the vector subspaces $K \subset L \subset M$ where dim$K=k$, dim$L=l$, and dim$M=m<n$. Take a sequence of $m$ linearly independent roots from $\Phi$ which uniquely determine $M$, such that the first $l$ of them uniquely determine $L$ and the first $k$ of them uniquely determine $K$. The permutation of the $m$ roots leads to different orientation of the cross-section with the $n$-polytope in $M$ and to different selection or permutation of the $k$ roots for $K$ and of the $l$ roots for $L$. Consider the permutations of $m$ roots that stabilize the $k$ roots in $K$ but do not stabilize the $l$ roots in $L$. As a result the cross-sections of the polytope in $L$ are changing but the cross-sections with $K$ remain the same. This provides with a good example of multiple polytopes in $l$ dimensions that have the same $k$ cross-section. For this same reason, a cross-section of a plane through the origin of a specifically oriented octahedron or cube both result in a hexagon. In the same way the observed polytope in the 3-subspace remains unchanged while many different orientations of the original polytope and cross-section polytopes can be observed in higher dimensions. \newline

Another result of even higher importance is related to hidden symmetries of the $n$-polytope that could not be obtained through its root system. Choosing three independent roots guarantees that the cross-section will pertain at least the symmetries resulting from these roots i.e. will be invariant under reflection through the hyperplanes orthogonal to the roots. However, the cross-section might have additional symmetries and could remain invariant under additional reflections that do not preserve the original $n$-polytope. We can easily illustrate this with an example. \newline

A cross-section of 3-cube with a 2-flat that intersects the middle points of six of its edges results in a hexagon. The symmetry that should formally be observed is the one of the $D_6$ group since all roots of the cube that are parallel to the 2-flat are long and generate $A_2$ root system. The angle between two roots parallel to the 2-flat is $2\pi/3$ and thus reflection along their respective hyperplanes (planes in this case) result in a 3-fold rotation that will preserve the cube. Another cross-section that is parallel to these two roots and intersects three vertices of the original cube gives an equilateral triangle where the extra symmetry is already not present. In the case of the initial cross-section, the roots of the $A_2$ root system all point to the vertices of the hexagon. As an even-sided polygon, the hexagon is also invariant under reflection defined by another root which points to the middle of its side.

\begin{figure}[h!]
  \centering
    \includegraphics[width=0.6\textwidth]{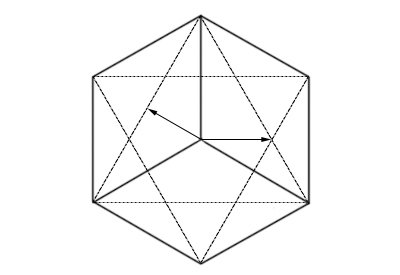}
  \caption{Hidden symmetry in the cube.}
\end{figure}

Denote the two simple roots of the $G_2$ root system $\alpha$ (pointing to a vertex and coinciding with a root from the required $A_2$ root system) and $\beta$ (pointing to the middle of a side). After observing the net of the cube, it can easily be concluded that all roots of the hexagon of type $\beta$ point to coordinates $(k/3,l/3),k=1,2,l=1,2$ on a given face of the cube, a total of four for each side. A general reflection of the cube along the $\beta$ vector fixes the orthogonal plane $H_{\beta}$ crossing two opposite vertices and the middle points of two edges. It does not preserve the cube and therefore does not belong to the symmetry group of the cube. As a result the $\beta$ vector cannot be a root of the $B_3$ root system and the $D_{12}$ group as the symmetry group of the hexagon is not a subgroup of the octahedral group $O_h$. Parabolic subgroups of the Coxeter group of a hypercube thus do not encompass all symmetries that might be observed in a cross-section with the hypercube and selecting roots from its root system only ensures the minimum symmetry of the particular cross-section. \newline

Symmetry of larger order can be obtained through cross-sections of higher dimensions in the same way. A central cross-section with a 4-cube that is parallel to three independent long roots generating a root system $A_3$ results in a specifically oriented octahedron (instead of a tetrahedron) which remains invariant under reflections along all roots in $A_3$. Hidden symmetries are also the reason why uniform polyhedra like the cuboctahedron cannot directly be obtained from intersections with an $n$-cube. Cuboctahedron can easily be obtained from a central cross-section with a 16-cell (4-orthoplex) that is parallel to three short roots of the $C_n$ root system and 16-cell can be obtained from a 4-dimensional central cross-section with the 5-cube. However, any further permutation of the sequence of four roots which are sufficient for a 16-cell cross-section will result in a $\binom{4}{3}$ possible 3-dimensional cross-sections none of which will take advantage of the additional symmetry of the 4-dimensional cross-section necessary for obtaining the cuboctahedron. \newline

A final conclusion then concerns the 3-dimensional cross-sections that can be obtained from an $n$-dimensional cubic tessellation. Rotating the tessellation will result in different cross-sections with its space-filling hypercubes and therefore different additional symmetry of the resulting space-filling polyhedra. Returning to the simple 3-dimensional example, the squares tiling can be obtained from a face-first intersection with the cubic honeycomb and the equilateral triangles from a vertex-first cross-section through three of the vertices of any of the cubes. The hexagon tiling however can be achieved only if the cubes are interpreted as polyhedra with additional vertices at the center of each original face with solid angles of $\pi$ and the original cubes are then arranged in a way so that the vertex of a cube touches the center of the original face of another and its edges remain parallel to the edges of the other. If the original cubic honeycomb is retained and at least one of the cubes is intersected at the centers of six of its edges, the resulting tessellation will be the trihexagonal tiling of the plane which consists of polygons with additional symmetry and polygons with the minimum necessary symmetry.

\newpage
\appendix
\appendixpage
\addappheadtotoc

\section{Visualizations}

\begin{figure}[h!]
  \centering
	$\begin{array}{cc}
		\includegraphics[scale=0.4]{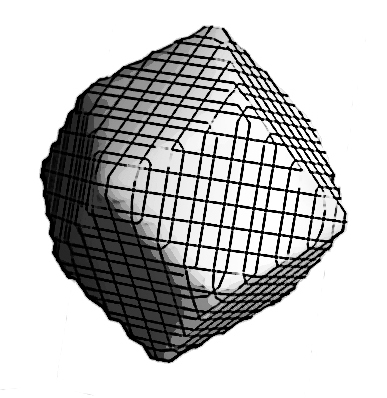} &
		\includegraphics[scale=0.4]{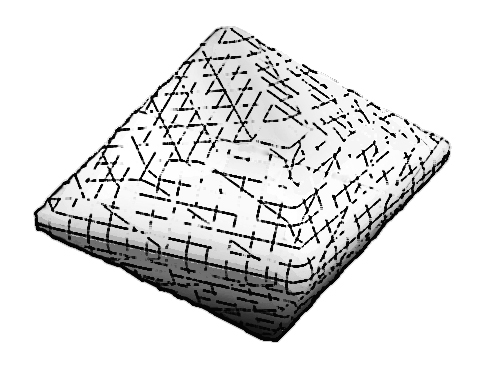}
	\end{array}$
  \caption{The convex polytope algorithm - rhombic dodecahedron obtained from 24-cell cross-section and hexagonal bipyramid obtained from 16-cell cross-section.}
\end{figure}

\begin{figure}[h!]
  \centering
	$\begin{array}{cc}
		\includegraphics[scale=0.4]{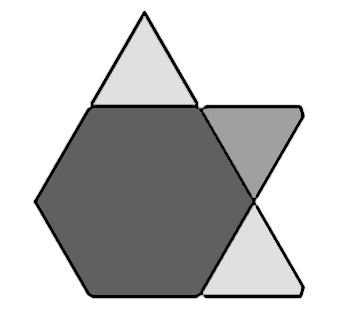} &
		\includegraphics[scale=0.4]{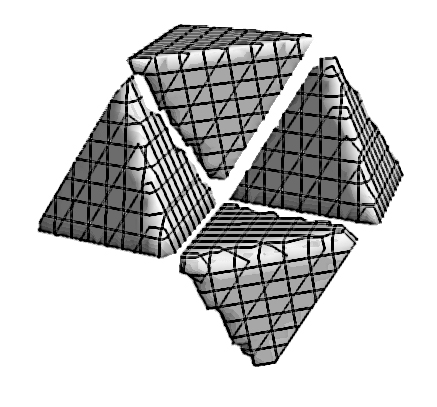}
	\end{array}$
  \caption{The hypercube tessellation algorithm - trihexagonal tiling and alternated cubic honeycomb.}
\end{figure}




\end{document}